\documentclass[10pt]{article} 
\usepackage{amsfonts}
\usepackage{amsmath}

\newcommand{\fsmpsp}{\left ( \Omega, \mathcal{F}, \{\mathcal{F}_t\}, P \right)}
\newcommand{\e}{\mathrm{E}}
\newcommand{\mcal}[1]{\mathcal{#1}}
\newcommand{\dx}{\mathrm{d}}
\newcommand{\ltwo}{{\bf{L}^2}}
\newcommand{\linf}{\bf{L}^{\infty}}
\newcommand{\lp}{\bf{L}^{p}}
\newcommand{\bb}[1]{\mathbb{#1}}
\newcommand{\norm}[1]{\left\lVert #1 \right\rVert}

\newcommand{\cxl}{\widetilde{\Lambda}}
\newcommand{\cxf}{\tilde{f}}
\newcommand{\hp}{{G}^{\perp}}
\newcommand{\cinth}{\mcal{Z}}
\newcommand{\dqdp}{\frac{\dx Q}{\dx P}}
\newcommand{\tw}{\widetilde{W}}
\newcommand{\ow}{\overline{W}}
\newcommand{\hcl}{\overline{\mcal{G}}}
\newcommand{\acl}{\overline{\mcal{A}}} 

\newcommand{\ltk}{B_k}
\newcommand{\Tau}{\tau}
\newcommand{\zq}{\mcal{Z}_Q}
\newcommand{\rhom}{\rho_{{G}}}

\newtheorem{lem}{Lemma}
\newtheorem{thm}{Theorem}

\newtheorem{asmp}{Assumption}

\newtheorem{notn}{Notation}
\newtheorem{prb}{Problem} 

\begin{document}

\title{Capital requirement for achieving acceptability}
\author{Soumik Pal\\ {\scriptsize{Department of Statistics, Columbia University}}\\ {\scriptsize{1255 Amsterdam Avenue, New York, NY 10027 (e-mail: soumik@stat.columbia.edu)}}}
\date{\today}

\maketitle

\begin{abstract}
Consider an agent who enters a financial market on day $t=0$ with an initial capital amount $x$. He invests this amount on stocks and the money market, and by day $t=\tau$, has generated a wealth $W$. He is given a convex class of probability measures or \emph{scenarios} and a real-valued \emph{floor} corresponding to each scenario. The agent faces the constraints that the expectation of $W$ under each scenario must not be less than the corresponding floor. We call $x$ \emph{acceptable} if one can start with $x$ and successfully generate $W$ satisfying these constraints.

The set of acceptable $x$ is a half-line in $\mathbb{R}$, unbounded from above. We show that under some regularity conditions on the set of scenarios and the floor function, the infimum of this set is given by the supremum of the floors over all scenarios under which $S$ is a martingale. 
\end{abstract}

\smallskip

\section{Introduction}\label{intro}

Consider an agent who trades on a time-interval $[0,\Tau]$ of finite length. The market offers finitely many assets; the agent invests an initial amount $x$, and by trading during the finite time horizon $[0,\Tau]$, ends up with an amount $W_{\Tau}$ at the end of the trading period. 
Now any pay-off one can generate at time $\Tau$ by trading in the market has a certain risk associated with it. We assume that one can measure this risk by a real number in a satisfactory manner. We shall call $W_{\Tau}$ \emph{acceptable}, if the risk of $W_{\Tau}$, which we denote by $\rho(W_{\Tau})$, is non-positive. In this paper we take up the problem of finding the minimum capital required to lead the agent, by careful trading, to an acceptable financial position at the end of the trading period.  

The notion of \emph{acceptability} of financial positions has been discussed in several recent papers. They are typically associated with defining a proper measure of risk. For example, in \cite{coherent}, Artzner, Delbaen, Eber \& Heath adopt an axiomatic set-up and introduce coherent measures of risk as real values associated with bounded random losses which satisfy certain desirable axioms. Interestingly, it turns out that the only coherent measures of risk are maximal expectations over a set of measures. Similarly in \cite{floor}, Carr, Geman and Madan extended the idea of acceptability to random variables representing potential gain from a derivative position. Acceptability of a such a random variable is characterized by a variety of measures, called \emph{valuation measures} and \emph{stress measures}. These authors define a random variable to be acceptable, if its expectation under each valuation measure is non-negative and its expectation under each stress measure is greater than or equal to a floor associated with the measure. F\"{o}llmer \& Schied introduce, in \cite{follmer}, the class of \emph{convex measures of risk}, which includes the above situations as special cases. We shall return to convex measures of risk in subsection~\ref{riskmeasure}. 

Finally, our own inspiration comes from a recently published paper \cite{shreve} in which the authors consider finitely many valuation and stress measures and associated floors, and characterize capital requirement from which one can trade to acceptability. In this paper, we extend their results to arbitrary families of such measures and floors, not necessarily finite or even countable.

The paper is divided as follows. Subsection~\ref{market} is devoted to ploughing through some technical grounds in order to present a precise statement of the problem. Curious readers can skim through such details, and have a look at subsection~\ref{statement} to get an idea about the statement of the problem. Our solution to this problem is through three main results in Section~\ref{original}, Theorems~\ref{mainthm},~\ref{minx}, and~\ref{finalthm}, with increasing ease of application at the cost of generality. Theorem~\ref{mainthm} gives a general analytical solution, which might be difficult to verify in practice. The probabilistic relevance becomes clear in Theorem~\ref{finalthm}, which exhibits how such a condition can be achieved by a proper choice of the underlying filtration. Finally in Section~\ref{examples} we use these results to compute capital requirements for efficiently hedging a claim, when we allow controlled shortfalls. 

\subsection{Acknowledgements}

My greatest debt of gratitude goes to Prof. Ioannis Karatzas who has patiently gone through the entire manuscript several times suggesting numerous corrections and improvements. This paper uses mathematical theory which I was unaware of to begin with. I sincerely thank Prof. Simeon Reich, Prof. Heinz Bauschke and Prof. Leonard Gross for suggesting me proper directions and references when I needed them. Thanks are also due to Prof. Peter Bank and Prof. Jaksa Cvitanic, who were kind enough to read through the final version of the manuscript and offered helpful comments.

\subsection{Description of the market}\label{market}

The market we consider has one risky asset and zero risk-free interest rate. These are simplifying assumptions, not difficult to avoid. But we adhere to them for notational simplicity. The price of our risky asset is assumed to be a real-valued (although only notational changes are required, in order to handle a vector-valued semimartingale) special semimartingale $S_t, \; 0\leq t \leq \Tau$, adapted to a suitable filtered probability space $\fsmpsp$. We assume that the filtration is right-continuous, $\mcal{F}_0$ contains all the $P$ null sets, and that $\mcal{F}_{\Tau}$ is the entire $\sigma$-algebra $\mcal{F}$. The semimartingale $S$ has the following Doob-Meyer decomposition
\begin{equation}\label{decomp}
S_t = S_0 + M_t + A_t, \;\; 0\leq t \leq \Tau,
\end{equation}
where the process $M$ is a local martingale and $A$ is a predictable process of finite variation. They are both assumed to be c\`{a}dl\`{a}g. Without loss of generality, at time zero, the initial price $S_0$ is assumed to be zero.

\noindent For any special semimartingale $X$ which can be decomposed as $X=N+V$, where $N$ is a local martingale and $V$ is a predictable finite variation process, one can define the $\mcal{H}^2$ norm of $X$ by $\norm{X}^2_{\mcal{H}}:= \e\left([N]_{\Tau} \right) + \e\left( \vert V \vert_{\Tau}^2 \right)$.
Here $[N]$ is the \emph{quadratic variation} of the local martingale $N$. The class of special semimartingales with a finite $\mcal{H}^2$ norm is a Banach space (see~\cite{protter}).

\begin{asmp}\label{finiteh2}
We shall assume that the $\mcal{H}^2$-norm of the semimartingale $S$ in~(\ref{decomp}) is finite, i.e., $\norm{S}_{\mcal{H}}:= \sqrt{\e\left([M]_{\tau} \right) + \e\left( \vert A \vert_{\Tau} ^2\right)} < \infty$.
\end{asmp}
Let $\Theta$ denote the collection of predictable processes $\pi$ such that 
\begin{equation}\label{h2normforpi}
\e\left(\int_0^{\Tau} \pi_u^2 \dx [M]_u\right) + \e\left(\int_0^{\Tau}\vert \pi_u\vert \dx \vert A \vert_u\right)^2 < \infty.
\end{equation}
Then for any predictable process $\pi \in \Theta$ and for any $0 \leq t \leq \Tau$, the stochastic integral of $\pi$ with respect to the process $S$ is well defined in the interval $[0,t]$ and will be denoted by $(\pi.S)_t := \int_0^t \pi_u \dx S_u$. The process $(\pi.S)_t$ is again a special semimartingale in the interval $[0,\tau]$ with a finite $\mcal{H}^2$ norm, whose square is given by~(\ref{h2normforpi}). See \cite{protter} for the proofs.

For $x \in \mathbb{R}$ and $\pi \in \Theta$, we call
\begin{equation}\label{wealth}
 W_u^{x,\pi} := x + (\pi.S)_u, \;\; 0\leq u \leq \Tau,
\end{equation}
the {\rm{wealth process}} at time $u$ starting with initial capital $x$ and generated by the trading strategy $\pi$. We shall make use of the following notation.

\begin{notn}\label{ltwot}
Let $\ltwo$ denote the space of all $\mcal{F}$-measurable random variables which are square-integrable under $P$. 
\end{notn}

Now for any $\pi \in \Theta$, it is again a standard fact that the random variable $(\pi.S)_{\Tau}$ is square integrable. Thus there is an obvious map from $\Theta$ into $\ltwo$ which carries $\pi$ to the stochastic integral $(\pi.S)_{\Tau}$. We consider the range of this map 
\begin{equation}\label{whatisg}
G := \left \{ X \in \ltwo \; \Big\vert \; X=\int_0^{\Tau} \pi_u \dx S_u, \;\;for \; some \;\; \pi\in \Theta \right \}.
\end{equation}
Clearly, $G$ is a subspace of the Hilbert space $\ltwo$. We shall denote by $\overline{G}$, the closure of $G$ in $\ltwo$. This closure will then be a Hilbert space in its own right. 

\subsection{Statement of the problem}\label{statement}

Let $\Delta$ be a collection of probability measures on $(\Omega, \mathcal{F})$, which are absolutely continuous with respect to $P$ and, let $\phi$ be a mapping from $\Delta$ into $\mathbb{R}$.

\begin{asmp}\label{deninl2}
Assume that $\Lambda\stackrel{\triangle}{=}\left\{\dx Q/\dx P\;\Big\vert\; Q\in\Delta \right\}$ is a subset of $\ltwo$.
\end{asmp}

\begin{prb}\label{problem1} Let $\Gamma$ be the subset of $\ltwo$ defined by
\begin{equation}\label{at}
\Gamma\stackrel{\triangle}{=}\left\{X\in\ltwo\;\;\Big\vert\;\;\e^Q\left( X \right) \geq \phi(Q), \;\; \forall\;Q \in \Delta\right\}.
\end{equation}
A real number $x$ will be called {\it acceptable} if
\begin{equation}\label{*t}
(x + G) \cap \Gamma \neq \emptyset,
\end{equation}
where the subspace $G$ is defined in~(\ref{whatisg}).
That is to say, $x$ is {\rm acceptable} if there exists a $\pi \in \Theta$ such that
\begin{equation}\label{*}
\e^Q\left( x + \int_0^{\Tau} \pi_u \dx S_u \right) \geq \phi(Q), \;\; \forall\;Q \in \Delta.
\end{equation}
It is immediate that the set of acceptable initial positions, which we shall denote by $\mcal{A}$,  is a half-line, unbounded from above. We shall be concerned with determining 
\begin{equation}\label{whatisinfa}
\inf\{x\in \mathbb{R} \vert \; x \in \mcal{A}\}.
\end{equation}
\end{prb}

\noindent{\sc remark.} Another important question is whether the set $\mcal{A}$ is closed or not. That is to say, whether the infimum in~(\ref{whatisinfa}) is attained. As we shall see in the beginning of Section~\ref{original}, our present set-up is deficient in answering the question. We shall get, however, a partial solution. 

\subsection{Measures of risk and Inf-convolutions.}\label{riskmeasure}

\emph{Convex measures of risk} have been introduced and discussed in detail in \cite{cxrisk}. The authors introduce axioms for convex measures of risk and show that coherent measures of risk, introduced in \cite{coherent}, are a subclass of the convex measures. The reader can find several examples of convex measures of risk in~\cite{follmer}.

For example, consider the space $\bf{L}^{\infty}$ of real-valued, $P$-essentially-bounded measurable functions defined on $({\Omega}, \mcal{F}, P)$ . Define the function $\rho$ by
\begin{equation}\label{rhox}
\rho \left(X \right) = \sup_{Q\in \Delta} \left( \e^Q\left[-X\right] + h\left( Q\right) \right ), \;\; X \in \bf{L}^{\infty},
\end{equation}
where $\Delta$ is as in the last section and $h:\Delta\rightarrow \mathbb{R}$. Then $\rho$ is a convex measure of risk, as discussed in page 172 of~\cite{follmer}.
Given the subspace $G$ of~(\ref{whatisg}), we can modify $\rho$ to obtain another measure of risk
\[
\rho_{{G}}(X) \stackrel{\triangle}{=} \inf_{H\in{G}}\rho(X-H),\; X\in\linf.
\]
This is an example of \emph{measures of risk in a financial market}, which is described in the discrete time setting in page 203 of \cite{follmer}. In general, these are special cases of \emph{inf-convolution} of risk measures developed in \cite{barrieu} and \cite{barrieu2}.

By Assumption~\ref{deninl2}, we can extend the domain of $\rho$ and $\rhom$ to the whole of $\ltwo$. Fix $\chi\in\ltwo$. We shall show that by a suitable choice of $\phi$, the value of the infimum in~(\ref{whatisinfa}) is equal to $\rhom(\chi)$. To see this, define ${\phi}(Q):=h(Q)-\e^Q(\chi)$. Then, by definition~(\ref{*t}), we get $\inf\{x\in\mathbb{R}\;\vert\;x\in\mcal{A}\}$
\begin{eqnarray}\label{nicole}
&=&  \inf\left\{ x\in\mathbb{R} \;\Big\vert\; \exists\; \xi\in G,\;\; \e^Q(x + \xi) \geq \phi(Q),\;\forall Q \in \Delta \right\}\nonumber\nonumber\\
&=& \inf\left\{ x \in\mathbb{R} \;\Big\vert\; \exists\; \xi\in G,\;\; \e^Q(x + \xi)   \geq h(Q)-\e^Q(\chi),\;\forall Q \in \Delta \right\}\nonumber\\
&=& \inf\left\{ x \in\mathbb{R} \;\Big\vert\; \exists\; \xi\in G,\;\; \e^Q(-(\chi + x + \xi)) + h(Q) \le 0,\;\forall Q \in \Delta \right\}\nonumber \\
&=& \inf\left\{ x\in\mathbb{R} \;\Big\vert\; \exists\; \xi\in G,\;\; \rho(\chi + x +\xi) \leq 0\right\}\nonumber\\
&=& \inf\left\{ x\in\mathbb{R} \;\Big\vert\; \inf_{H\in G}\rho(\chi+x-H) \le 0\right\}\nonumber\\
&=& \inf\left\{ x\in\mathbb{R} \;\Big\vert\; \rhom(\chi + x) \le 0\right\}=\rhom(\chi).
\end{eqnarray}
The fifth equality above requires proper assumption on the regularity of $\rho$ and the last one is due to \cite{follmer}, page 155, eqn.(4.5).

\section{A general Hilbert space problem}\label{hilbertspace} 

Let $\mcal{H}$ be a Hilbert space with an inner product denoted by $\langle . \;, .\; \rangle$ and the norm by $\norm{.}$. Suppose we are given a set $\Lambda \subseteq \mcal{H}$, a mapping
$f:\Lambda \rightarrow \mathbb{R}$, and a closed subspace $\mcal{G}\subseteq \mcal{H}$. For any given real number $x$, we want to find necessary and sufficient conditions for the existence of an element $z^*\in \mcal{G}$ such that
\begin{equation}\label{**}
\langle z^*, y \rangle \geq f(y) - x , \;\; \forall y\: \in \Lambda.
\end{equation}
The rest of this section is devoted to solving this problem.

Let $\cxl$ denote the convex hull of $\Lambda$. One can extend the mapping $f$ from $\Lambda$ to $\cxl$ by defining a new mapping $\cxf:\cxl \rightarrow \mathbb{R}$, given by
\begin{equation}\label{whatiscxf}
\cxf(y) := \sup \left \{ \sum^n_{i=1} \lambda_i f(z_i) \right \},\; \;y \in \cxl.
\end{equation}
Here the supremum is taken over all choices of $n \in \mathbb{N}$ and $\lambda_1 \geq 0, \ldots, \lambda_n \geq 0$ with $\sum_i \lambda_i =1$) and $z_{1}\in \Lambda, \ldots, z_{n} \in \Lambda$ which satisfy $y = \sum_{i=1}^n \lambda_i z_{i}$. 

Let us observe that if $z^*$ is a solution for~(\ref{**}), then $z^*$  also solves a more general class of inequalities. In fact, by the linearity of inner products, it follows from~(\ref{**}) that if $y\in\cxl$ can be written as as a convex combination of some $\{z_1, z_2,\ldots,z_n\} \subseteq \Lambda$, i.e. $y=\sum\lambda_iz_i$, then $\langle z^*, y \rangle = \sum \lambda_i.\langle z^*, z_i \rangle \geq \sum\lambda_if(z_i) - x$. Thus, we can appeal to the definition of $\cxf$ in~(\ref{whatiscxf}) to obtain
\begin{equation}\label{cx**}
\langle z^*, y \rangle \geq \cxf(y) - x, \;\; \forall y \in \cxl.
\end{equation}

Let $T(y)$ for any $y\in \mcal{H}$ denote the unique orthogonal projection of $y$ on $\mcal{G}$. In particular, we have
\begin{equation}\label{ty}
\langle z , y \rangle = \left\langle z , T(y) \right\rangle, \;\; \forall z \in \mcal{G}.
\end{equation}

\begin{thm}\label{maininhs}
For any $x\in \mathbb{R}$, a necessary and sufficient condition for the existence of $z^* \in \mcal{G}$ satisfying the inequalities in~(\ref{**}), is the existence of a constant $M \geq 0$ such that
\begin{equation}\label{maineqn}
M \norm{T(y)} \geq \cxf(y) - x, \;\; \forall y \in \cxl.
\end{equation}
\end{thm}

\noindent{\sc{proof.}} To see the necessity of condition~(\ref{maineqn}), just apply the Cauchy-Schwarz inequality to~(\ref{cx**}) to get
\begin{equation}\label{karat}
\cxf(y) - x \leq \langle z^*, y \rangle = \left\langle z^*, T(y) \right\rangle \leq \norm{z^*}\norm{T(y)}, \;\; \forall y \in \cxl.
\end{equation}
Setting $M := \norm{z^*}$ we have established condition~(\ref{maineqn}). 

Proving the sufficiency is more subtle. We start with the assumption that~(\ref{maineqn}) holds for some $x\in \mathbb{R}$ and some real constant $M\geq 0$. To simplify notation, let us define a new mapping $b:\cxl \rightarrow \mathbb{R}$ by
\[
b(y) := \cxf(y) - x, \;\; y \in \cxl. 
\]
Condition ~(\ref{maineqn}) then reads
\begin{equation}\label{condforb}
M\norm{T(y)} \geq b(y), \;\; \forall \; y \in \cxl.
\end{equation}
We shall establish~(\ref{**}) by showing that there exists $z^*\in \mcal{H}$ such that
\begin{equation}\label{reformof13}
\langle z^*, y \rangle \geq b(y), \;\; \forall\; y \in \cxl.
\end{equation}

\noindent$\bullet$ We shall first show that for any given finite subset $\{y_1,\ldots, y_n \} \subseteq \cxl$, there is a $z^* \in \mcal{G}$ such that $\norm{z^*} \leq M$ and $z^*$ satisfies
\begin{equation}\label{lem1}
\langle z^*, y_k \rangle \geq b(y_k), \;\; \forall \; 1 \leq k \leq n.
\end{equation}

We shall argue this by contradiction. Suppose that no such $z^*$ exists. Consider the set
\[
\mcal{S} := \left \{ \left ( \; \langle z, y_1 \rangle, \ldots, \langle z, y_n \rangle\; \right ) \Big\vert \; z \in \mcal{G},\; \norm{z} \leq M \right \}
\]
which is compact and convex in $\mathbb{R}^n$. Here and throughout, $\bb{R}^{n}_{+}$ will refer to the subset of points in $\bb{R}^n$ which have all co-ordinates non-negative. Let $\mcal{S}^-$ be the set all points $(a_1, a_2, \ldots, a_n)$ which can be represented as
\[
a_k = \langle z, y_k \rangle - r_k, \;\; 1\leq k \leq n,
\]
for some $r_k \geq 0$ and some $z \in \mcal{G}$ such that $\norm{z} \leq M$. For notational simplicity, let us denote $b_k:= b(y_k), \;\; 1\leq k \leq n$. Since we have assumed that no solution to ~(\ref{lem1}) exists, we have 
\begin{equation}\label{closed}
(b_1, \ldots, b_n) \notin \mcal{S}^-.
\end{equation}

But, by the Separating Hyperplane Theorem,~(\ref{closed}) implies that there exists a vector $\lambda=(\lambda_1, \ldots, \lambda_n) \in \bb{R}^n$, $\lambda\neq 0$, such that for all $a_1 \geq 0, \ldots, a_n \geq 0$ we have
\begin{equation}\label{sep}
\sum \lambda_i b_i > \sum \lambda_i \langle z, y_i \rangle - \sum \lambda_i a_i, \;\; \forall z \in \mcal{G}, \; \norm{z} \leq M.
\end{equation}
For any $i$, let $a_i$ tend to infinity to see that $\lambda_i \geq 0$. We can thus normalize $\lambda$ to have $\sum_{i=1}^n \lambda_i = 1$. Taking $a_1=\ldots= a_n=0$ in ~(\ref{sep}), we obtain
\begin{equation}\label{contraeq}
\sum \lambda_i b_i > \sum \lambda_i \langle z, y_i \rangle
= \left \langle z, \sum \lambda_i y_i \right \rangle = \left \langle z, y_c \right \rangle,
\end{equation}
where $y_c := \sum \lambda_i y_i \in \cxl$. 

Note that the function $\cxf$ in~(\ref{whatiscxf}) is concave. This has been proved in \cite{rock}, page 37, Theorem 5.6, and the example following right after its proof. Thus the function $b(.)=\cxf(.) - x$ is also concave. Combining with ~(\ref{condforb}), and using the concavity of $b$, we get
\begin{equation}\label{contramain}
M\norm{T(y_c)} \geq b(y_c) = b\left(\sum \lambda_i y_i \right) \geq \sum \lambda_i b_i > \langle z, y_c \rangle = \langle z, T(y_c) \rangle,
\end{equation}
for every $z\in\mcal{G}$ with $\norm{z} \leq M$. If $\norm{T(y_c)} = 0 $, this leads to $0 \geq \sum \lambda_i b_i > 0$, which is a contradiction; whereas if $\norm{T(y_c)} \neq 0$, note that $z= {M}.T(y_c)/\norm{T(y_c)}$ is an element of the subspace $\mcal{G}$ with $\norm{z} \leq M$ which, when plugged into inequality~(\ref{contramain}), gives 
\[
M \norm{T(y_c)} \geq \sum \lambda_i b_i > M \norm{T(y_c)},
\]
again a contradiction. We have thus proved ~(\ref{lem1}).
\medskip

\noindent$\bullet$ In general, let us define for any $y \in \cxl$, the following subset of $\mcal{G}$: 
\[
\Pi_{y} := \left \{ z \in \mcal{G}\; \Big\vert \; \norm{z} \leq M, \; \langle z, y \rangle \geq b(y) \right \}.
\]
Then there is a solution to ~(\ref{reformof13}) if we can show that 
\begin{equation}\label{nonemp}
\cap_{y \in \cxl} \Pi_{y} \neq \emptyset.
\end{equation}
Now each $\Pi_{y}$ is a closed subset of the $M$-ball of $\mcal{G}$, a set which is compact under the weak topology. This follows from the Banach-Alaoglu Theorem and the fact that a Hilbert space is its own dual (see \cite{rudin}, pages 68, 94). Thus $\Pi_{y}$ is a weak-compact subset of $\mcal{G}$. Hence, by the finite intersection property,~(\ref{nonemp}) holds if and only if for any finite collection $\{ y_1, \ldots, y_n \} \subseteq \cxl$, we have $\cap_{1\leq i \leq n} \Pi_{y_i} \neq \emptyset$.
That is to say,~(\ref{nonemp}) holds if and only if for any finite collection $\{ y_1, \ldots, y_n \} \subseteq \cxl$ we can find an element $z^*\in \mcal{H}$ such that $\norm{z^*} \leq M$ and $\langle z^*, y_k \rangle \geq b(y_k), \;\;  k=1,\ldots, n$.
But this is what we have shown in~(\ref{lem1}). This proves the theorem. \hfill $\Box$
\medskip

Our previous result does not hold when the subspace $\mcal{G}$ is not a closed subspace of the Hilbert space $\mcal{H}$. However what we shall show now is that not much is lost if we consider $\hcl$, the closure of $\mcal{G}$ instead of $\mcal{G}$ itself.

Let us denote by $\acl$ the set of all real numbers $x$ for which the inequalities in~(\ref{**}) have a solution for some $z^* \in \hcl$, and reserve the notation $\mcal{A}$ for that subset of $\acl$ for which the solution $z^*$ is actually an element of $\cal{G}$. We shall now show that when $\Lambda$ is bounded in norm, $\mcal{A}$ is a dense subset of $\acl$. However, since both $\mcal{A}$ and $\acl$ are half-lines, this is actually equivalent to proving what we shall need most, i.e.,
\begin{equation}\label{openh}
\inf \mcal{A} = \inf \acl.
\end{equation}
This is achieved by the following lemma.

\begin{lem}\label{lambdabdd}
If the set $\Lambda$ is bounded in norm and if $\;x \in \acl$, then 
$(x + \epsilon)\in \mcal{A}$ for any positive $\epsilon$. 
\end{lem}

\noindent{\sc{proof.}} Fix $x \in \acl$ and an $\epsilon > 0$. By the definition of $\acl$, there exists $z \in \hcl$ such that
\[
\langle z, y \rangle \geq f(y) - x, \;\; \forall y \in \Lambda.
\] 
Now since $\Lambda$ is bounded in norm and $\mcal{G}$ is dense in $\hcl$, there is an element $z^* \in \mcal{G}$ such that
\[
\sup_{y\in\Lambda}\vert \langle z^*, y \rangle - \langle z, y \rangle \vert \leq \norm{z^* - z}.\sup_{y\in\Lambda}\norm{y}\leq \epsilon.
\]
Hence, we get $\langle z^*, y \rangle \geq \langle z, y \rangle - \epsilon \geq f(y) - (x + \epsilon), \;\; \forall y \in \Lambda$.
Since $z^* \in \mcal{G}$, this shows that $(x+\epsilon) \in \mcal{A}$, and proves the lemma. \hfill $\Box$
\medskip

\noindent{\sc{remark.}} Equation~(\ref{openh}) does not hold in full generality, although we always have $\inf \mcal{A} \geq \inf \acl$,
since $\mcal{A}\subseteq \acl$. We shall return to discuss this point again in the next section.

\section{Main results}\label{original} 

We shall now translate the results of the last subsection in order to solve~(\ref{whatisinfa}). Consider the Hilbert space $\ltwo$ and the subspace $G$ of stochastic integrals defined in~(\ref{whatisg}). Let $\overline{G}$ denote the closure of $G$ in $\ltwo$. Recall the statement of the problem in subsection~\ref{statement}, and as in the setting of the last section, define:
\begin{equation}
\mcal{H}= \ltwo,\;\;\mcal{G} = \overline{G},\;\; \Lambda = \left\{\dx Q/\dx P\;\Big\vert \; Q\in\Delta\right\}. 
\end{equation}
That the set $\Lambda$ is a subset of $\mcal{H}$ is a consequence of Assumption~\ref{deninl2}. As before, $\cxl$ will denote the convex hull of $\Lambda$. Note that there is a one-to-one correspondence between the elements in $\Lambda$ and the probability measures in $\Delta$. Define the function $f:\Delta\rightarrow\mathbb{R}$ by
\begin{equation}\label{defnf}
f(X) = \phi(Q), \;\; \mathrm{for}\; X=\dx Q/\dx P,\; Q\in\Lambda.
\end{equation}
Define $\cxf$ on $\cxl$ in the same way as in~(\ref{whatiscxf}). The notation for $\norm{.}$, from now on, is strictly reserved for the $\ltwo$ norm.

Clearly with this set-up, for any $X \in \overline{G}$ and any measure $Q$ such that $\dx Q/ \dx P \in \cxl$, one has $\e^Q(X) = \langle X, \dx Q / \dx P \rangle$. This association makes evident the relation between solving inequalities~(\ref{*}) and~(\ref{**}). In fact, if $G$ is a closed subspace of $\ltwo$, solving for~(\ref{*}) is exactly the same as solving for~(\ref{**}). Problems arise when $G$ is not closed; for then the solution obtained in ~(\ref{**}) might be an element strictly in the closure of $G$. This problem is easy to deal with when $\Lambda$ is bounded in norm, since our object of interest, $\inf \mcal{A}$, remains the same whether we consider $G$ or $\overline{G}$, as shown by Lemma~\ref{lambdabdd} at the end of the last section. 

In general, however, we cannot expect that the wealth process $\int_0^t \pi_u\dx S_u$ which satisfies inequalities~(\ref{*}) will have finite $\mcal{H}^2$ norms. A good analogy will be to think of situations where the optimal wealth process is a strict local martingale instead of being a true martingale. Our subspace $G$ only allows terminal wealth from a wealth process which has finite $\mcal{H}^2$-norm, and this is usually a strong requirement. Thus it seems necessary that we reformulate Problem~\ref{problem1} by allowing solutions which belong to $\overline{G}$ rather than $G$ itself. We now restate Problem~\ref{problem1} in the following way:

\begin{prb}\label{problem2}
A real number $x$ will be called {\rm weakly acceptable} if
\begin{equation}\label{***}
(x + \overline{G}) \cap \Gamma \neq \emptyset,
\end{equation}
where the subspace $G$ is defined in~(\ref{whatisg}) and $\Gamma$ is defined in~(\ref{at}).

We shall denote the set of weakly acceptable numbers by $\acl$, which is still a half-line not bounded above. As before, we shall be concerned with determining $\inf\acl$.
\end{prb}

The operator $T$ will denote projection onto the subspace $\overline{G}$. That is, for any $X \in \ltwo$, one has the following decomposition:
\begin{equation}\label{zalpha}
X = T(X) + \left[ I - T \right](X),
\end{equation}
where $\left[ I - T \right](X)$ is orthogonal to every element in $\overline{G}$. The next theorem is a restatement of Theorem~\ref{maininhs}.

\begin{thm}\label{mainthm}
Under Assumptions ~\ref{finiteh2} and ~\ref{deninl2}, a real number $x$ is weakly acceptable (in the sense of Problem~\ref{problem2}) if and only if there exists a non-negative real constant $M$ such that
\begin{equation}\label{eqninmainthm}
M\norm{T(X)} \geq \cxf(X) - x, \;\;\forall \;X\in \cxl.
\end{equation}
\end{thm}

The probabilistic interpretation of $T(X)$ will be clear in the next lemma.
\begin{lem}\label{zeroproj}
For any $X\in\ltwo$, consider the process $X_t \stackrel{\triangle}{=} \e\left(X\;\Big\vert\;\mcal{F}_t\right),\;\; 0\leq t\leq \Tau$. Then $T(X)=0$ implies $\left\{X_t.S_t, \mcal{F}_t\right\}_{0\leq t\leq \Tau}$ is a martingale.
\end{lem}

\noindent{\sc{proof.}} For any stopping time $\sigma$ taking values in $[0,\tau]$, consider the process $\pi_u := 1_{\{\sigma \geq u\}}, \;\;0\leq u\leq \tau$. Since $T(X)=0$, we have $\e\left(X\int_0^{\Tau}\pi_u\dx S_u\right)=\e\left(XS_{\sigma}\right)=0$. Thus $\e\left(XS_{\sigma}\right)$ is zero for all stopping times $\sigma$. By taking conditional expectation with respect to $\mcal{F}_{\sigma}$, we have $\e(X_{\sigma}S_{\sigma})=0$ for all stopping times $\sigma$. This proves the lemma. \hfill$\Box$

\begin{lem}\label{zeroismgle}
For any $Q\in\Delta$, let $\zq=\dx Q/\dx P$ denote the Radon-Nikodym derivative of $Q$ with respect to $P$. Then $T(\zq) = 0$ if and only if the process $S$ is a $Q$-martingale on the interval $[0,\tau]$.
\end{lem}
\medskip

\noindent{\sc{proof.}} The \emph{only if} part follows from the last lemma via what is commonly know as the {\rm Bayes rule}. See, for example, \cite{karatshreve}, page 193.

For the \emph{if} part, start with a measure $Q$ such that the process $S_t$ is a martingale on the interval $[0,\Tau]$. One can show by an application of the Burkholder-Davis-Gundy inequality (a proof can be found in \cite{shreve}, Proposition 1) that under Assumption ~\ref{finiteh2}, for any $\pi \in \Theta$, the process $\int_0^t \pi_u \dx S_u$ is a martingale under $Q$. Thus 
\[
\e\left[ \frac{\dx Q}{\dx P} \int_0^{\tau} \pi_u\dx S_u\right]=\e^Q\left (\int_0^{\Tau} \pi_u \dx S_u \right ) = 0, \;\; \forall \pi \in \Theta.
\] 
This shows that $\zq$ is orthogonal to $G$ and hence $T(\zq)=0$.\hfill $\Box$

Hereafter a \emph{martingale measure} will refer to a probability measure $Q$ under which the process $\{S_u\}$ is a martingale in the interval $[0,\Tau]$. Let $\hp$ denote the orthogonal complement of the subspace $G$, defined in~(\ref{whatisg}). Then, by what we have just proved in the last lemma, the set 
\begin{equation}\label{cinthperp}
\cinth:={\cxl}\cap G^{\perp} = \left\{X\in \cxl \;\Big\vert \;T(X)=0\right\}
\end{equation}
is the set of probability densities corresponding to the martingale measures in the convex hull of $\Delta$. 

From now on we shall also assume the following.
\begin{asmp}
The function $\cxf$ is contiuous with respect to the $\ltwo$ norm on $\cxl$. 
\end{asmp}
One can then extend $\cxf$ continuously to the closure of $\Lambda$, which we shall, by an abuse of notation, continue to denote by $\cxl$. Our next theorem considers the case when $\cinth=\emptyset$, while the other case is taken up in Theorem~\ref{minx}.

\begin{thm}\label{nozero}
Suppose $\cinth = \emptyset$. If we have $\sup_{Q \in \Delta} \phi(Q) < \infty$, or equivalently,
\begin{equation}\label{finitef}
\sup_{X \in \Lambda} f(X) < \infty,
\end{equation}
then for $\inf\acl=-\infty$, where $\acl$ is defined in Problem~\ref{problem2}.
\end{thm}

In other words, under the condition~(\ref{finitef}), the non-existence of martingale measures in the closed convex hull of the set of scenarios, $\Delta$, implies that every $x\in\mathbb{R}$ is an weakly acceptable initial position.
\smallskip

\noindent{\sc{proof.}} The set $T(\cxl)$, the image of $\cxl$ under the orthogonal projection mapping $T$, is closed and convex. Since $\cinth = \emptyset$, we have $0 \notin T(\cxl)$. Thus a basic fact from Hilbert space theory states that there is an element in $T(\cxl)$ which is of minimum positive norm. That is, there is an element $X^* \in \cxl$ such that $0 < \norm{T(X^*)} = \inf_{X\in \cxl}\norm{T(X)}$. 

Note that~(\ref{finitef}) implies $\sup_{X \in \cxl} \cxf(X) < \infty$. One can then define $ K = \max(\sup_{X \in \cxl} \cxf(X) - x, 0)$, and consider $M = K/\norm{T(X^*)}$ to see that
\[
M\norm{T(X)}\geq M\norm{T(X^*)} = K \geq \cxf(X) - x, \;\; \forall \; X \in \cxl.
\]
This shows that~(\ref{eqninmainthm}) is satisfied, and proves the theorem.\hfill $\Box$
\smallskip

For any $X \in \ltwo$ and for any $1 \leq p \leq 2$, let us denote the $\lp$ norm of $X$ by $\norm{X}_p$, i.e., 
\[
\norm{X}_p := \left[\e(\vert X\vert^p) \right]^{1/p}.
\]
Since $X$ has finite second moment and we are on a probability space, an application of H\"{o}lder's inequality shows that $\norm{X}_p$ is finite for any $1\leq p \leq 2$. We also define the $\lp$-distance between a point $X\in\lp$ and a non-empty subset $\Pi\subseteq \lp$ by
\[
d_p(X,\Pi) \stackrel{\triangle}{=} \inf_{Y \in \Pi} \norm{X-Y}_p. 
\]
Again, the distance is well-defined and finite for any $1 \leq p \leq 2$. 

\begin{thm}\label{minx} Suppose the following assumptions are satisfied:

\begin{enumerate}
\item $\cinth \neq \emptyset$.

\item There exists a constant $L > 0$ and some $p\in (1,2]$ such that 
\begin{equation}\label{lplip}
\vert \cxf(X) - \cxf(Y)\vert \leq L\norm{X - Y}_p \;\; \forall \; X, Y \in \cxl.
\end{equation}

\item For any sequence $\{X_n\}\subseteq \cxl$ such that $\lim_{n\rightarrow 0}\norm{T(X_n)}= 0$, we also have (at least through a subsequence)
\begin{equation}\label{pcont}
\lim_{n\rightarrow \infty} d_p(X_n, \cinth ) = 0.
\end{equation}
\end{enumerate} 
Then we can conclude that 
\begin{equation}\label{firstinf}
\inf \acl=\sup_{Y \in \cinth} \cxf(Y). 
\end{equation}
Here $\acl$ is the set of weakly acceptable initial positions described in Problem~\ref{problem2}. 

\end{thm}

Our proof will be achieved by the following two lemmas. The first one needs the concept of \emph{nearest point projections} in uniformly convex (or uniformly rotund) Banach spaces, e.g. the $\lp$ spaces, $p\in(1,\infty)$. This can be found in~\cite{megg}, page 427, Example 5.1.4. We can then use corollary 5.1.19 on page 435 of~\cite{megg}, to see that given any closed, convex subset $\Pi$ and any element $X$, both in $\lp$ for some $1< p <\infty$, there is an element $S_{\Pi}(X) \in \Pi$ such that  
\begin{equation}\label{spi}
\norm{X - S_{\Pi}(X)}_p = \inf_{Y\in \Pi}\norm{X - Y}_p = d_p(X, \Pi). 
\end{equation}
Additionally, the operator $S_{\Pi}$ is \emph{sunny}, i.e., satisfies (see \cite{npp}, page $17$) 
\begin{equation}\label{scaling}
S_{\Pi}\big(\alpha X + (1-\alpha) S_{\Pi}(X)\big)= S_{\Pi}(X), \;\;\forall\; \alpha \geq 0.
\end{equation}
In what follows, we shall consider $\Pi$ to be the closure of $\cinth$ in $\lp$. Since $\cinth$ is $\lp$ dense in $\Pi$, it follows that any real function, uniformly continuous on $\cinth$ with respect to the $\lp$ metric, can be extended uniquely on $\Pi$. By our second assumption in Theorem~\ref{minx}, the function $\cxf :\cxl\rightarrow\mathbb{R}$ is uniformly continuous with respect to the $\lp$ and hence can be extended to elements of $\Pi$. 

The proofs of the following lemmas are done in the appendix.

\begin{lem}\label{Aisbdd}
Under the assumptions and notation of Theorem~\ref{minx}, there exists a constant $M_1 \in [0,\infty)$ such that
\begin{equation}\label{plipt}
d_p(X,\cinth) \leq M_1 \norm{T(X)}, \;\; \forall\; X \in \cxl.
\end{equation}
\end{lem}

\begin{lem}\label{Bisbdd}
For a given $z\in \mathbb{R}$, suppose there exists a constant $M_2 \in [0,\infty)$ (may depend on $z$), such that  
\begin{equation}\label{blipt}
\cxf(S_{\Pi}(X)) - z \leq M_2\norm{T(X)}, \;\;\forall\; X \in \cxl;
\end{equation}
then $z \geq \sup_{X\in\cinth}\cxf(X)$. Conversely, for any $z\geq \sup_{X\in\cinth}\cxf(X)$, clearly~(\ref{blipt}) holds with $M_2=0$.
\end{lem}

\noindent{\sc{proof of theorem~\ref{minx}.}} Choose $x\in\mathbb{R}$. For any $X\in \cxl$, one has the decomposition
\begin{equation}\label{decbbnd}
\cxf(X) - x = \cxf(X)- \cxf\left(S_{\Pi}(X)\right) + \cxf\left(S_{\Pi}(X)\right) - x.
\end{equation}
By Lemma ~\ref{Aisbdd}, there is a $M_1\in [0,\infty)$ such that
\[
d_p(X,\cinth)=\norm{X - S_{\Pi}(X)}_p \leq M_1\norm{T(X)},
\]
and thus, by assumption (2) in Theorem~\ref{minx}, we obtain
\begin{equation}\label{firstisbdd}
\vert \cxf(X)- \cxf\left(S_{\Pi}(X)\right)\vert \leq L\norm{X - S_{\Pi}(X)}_p \leq L.M_1\norm{T(X)}.
\end{equation}
Plugging in the above inequality in~(\ref{decbbnd}), we see that~(\ref{eqninmainthm}) holds, for some $M\geq 0$, if and only if there exists a constant $M_2$ for which
\[
\cxf\left(S_{\Pi}(X)\right) - x \leq M_2\norm{T(X)}, \;\;\forall X\in \cxl.
\]
But by Lemma~\ref{Bisbdd}, this can happen if and only if $x \geq \sup_{Y \in \cinth} \cxf(Y)$. This shows that $\inf\acl = \sup_{Y \in \cinth} \cxf(Y)$ and proves Theorem~\ref{minx}. \hfill $\Box$
\medskip

Condition~(\ref{pcont}) of Theorem~\ref{minx} will, in general, not be easy to verify. However, our next result displays an interesting link between the geometric and probabilistic aspects of the problem. It shows that with an appropriate assumption on the underlying filtration of the stock price process, we can make~(\ref{pcont}) automatic.  

\begin{thm}\label{finalthm}
Let $\Delta$ of subsection~\ref{statement} be the set of all probability measures $Q$ on $(\Omega, \mcal{F})$, such that $Q \ll P$ and $\norm{{\dx Q}/{\dx  P}} \leq \mcal{K}$, for some given constant $\mcal{K}\in(0,\infty)$. As before, $\cxl$ will denote the collection of Radon-Nikodym derivatives of the measures in $\Delta$, i.e.,
\begin{equation}\label{describecxl}
\cxl \stackrel{\triangle}{=} \left\{ X\in\ltwo \;\Big\vert\; X\geq 0\; \mathrm{a.s.}\;P, \;\e(X)=1\; {\mathrm {and}} \norm{X}\leq \mcal{K} \right\}.
\end{equation}
Assume that 
\begin{enumerate}
\item there exists an element $M^*$ of $\cinth$, defined in~(\ref{cinthperp}), with $\norm{M^*} < \mcal{K}$;
\item all martingales of the filtration $\{\mcal{F}_t\}$ have continuous versions; and 
\item the mapping $\cxf$ satisfies~(\ref{lplip}). 
\end{enumerate}

Then, we have $\inf \mcal{A} = \sup_{X \in \cinth} \cxf(X)$. 
\end{thm}

To prove this theorem, we shall show that the assumptions in Theorem~\ref{minx} are satisfied. In fact, we only need to show that~(\ref{pcont}) holds, since the other assumptions are already assumed to be true. The proof of the following lemma is in the appendix.

\begin{lem}\label{approxmart} Let $\{Y_n\}$ be a sequence in $\cxl$ such that $\lim_{n\rightarrow\infty}T(Y_n)=0$. Then there exists a sequence $\{L_n\} \subseteq G^{\perp}$, with $P(L_n \geq 0)=1$ and $\e(L_n)=1$, such that:

\begin{equation}\label{stepone1}
\lim_{n\rightarrow \infty}\norm{Y_n - L_n}_p = 0
\end{equation}
and
\begin{equation}\label{stepone2} 
\limsup_{n\rightarrow \infty}\norm{L_n}_2 \leq  \mcal{K}. 
\end{equation}

\end{lem}

\noindent{\sc{proof of theorem~\ref{finalthm}.}} Consider $M^*$ as in assumption 1 of Theorem~\ref{finalthm} and the sequence $\{L_n\}$ from Lemma~\ref{approxmart}. For any $\alpha \in (0,1)$, define the sequence $W_n:=\alpha L_n + (1 - \alpha)M^*$. Then, by the triangle inequality, we have 
\[
\limsup_{n\rightarrow\infty}\norm{W_n}\leq \alpha\limsup_{n\rightarrow \infty}\norm{L_n} + (1-\alpha)\norm{M^*}.
\]
Thus, from~(\ref{stepone2}), we get that $\limsup_{n\rightarrow\infty}\norm{W_n} < \mcal{K}$. In other words, there is a large enough $N$ such that $\norm{W_n} < \mcal{K}$ for all $n > N$. Now, by Lemma~\ref{approxmart}, each $P(L_n\ge 0)=1$ and integrates to one. Also, since $M^*$ is a probability density, $P(M^*\ge 0)=1$ and $\e(M^*)=1$. Thus we also have $P(W_n \geq 0)=1$ and $\e(W_n)=1$, and thus from~(\ref{describecxl}), $W_n \in \cxl, \;\; \forall\; n > N$. But, again by Lemma~\ref{approxmart}, each $L_n$ belongs to $G^{\perp}$. Since $M^*$ also belongs to $G^{\perp}$, we conclude 
\[
W_n \in \cxl \cap G^{\perp} = \cinth,\;\;\forall n > N.
\]
Clearly then
\begin{eqnarray}\label{steptwo}
\limsup_{n\rightarrow \infty} d_p(Y_n, \cinth) &\leq& \limsup\norm{Y_n- W_n}_p \nonumber\\
 &\leq& \limsup\norm{Y_n- \alpha L_n - (1 - \alpha)M^*}_p \nonumber\\
 &\leq& \limsup\norm{Y_n - L_n}_p + (1-\alpha)\limsup\norm{L_n - M^*}_p\nonumber\\
 &=& 0 + (1 - \alpha)\left(\limsup\norm{L_n}_p + \norm{M^*}_p\right)\nonumber\\
&\le& (1 - \alpha)\left(\limsup\norm{L_n} + \norm{M^*}\right)\le (1 - \alpha)2\mcal{K}.
\end{eqnarray} 
The final inequality is due to~(\ref{stepone2}) while the one right before it follows from H\"{o}lder's inequality: for any random variable $Z$, we have  
\begin{equation}\label{holderyield}
\norm{Z}_p \leq \norm{Z}_2 = \norm{Z}, \;\;\forall\; 1 < p < 2.
\end{equation}  
Take $\alpha \uparrow 1$ in the above inequality to conclude that $\limsup d_p(Y_n, \cinth) = 0$ which shows~(\ref{pcont}) holds and the proof of Theorem~\ref{finalthm} is complete.\hfill$\Box$

\section{Examples}\label{examples} We solve a prototypical example of determining the sellers' price of an option in an incomplete market. Due to incompleteness of the market, a typical contingent claim will not admit a perfect hedge. In the following example, we show that instead of taking the conservative apporach of superhedging, an investor can allow some controlled risk of shortfall, and then compute the necessary initial capital for an {\it efficient} hedge.

\noindent {\bf Example 1.} Consider a market with two stocks whose price processes $S$ and $S^{'}$ are driven by a two dimensional Brownian motion till a finite terminal time $\Tau$. For simplicity, the rate of interest, the mean rate of return, and the rate of dividend are kept at zero. The price process $S$ of stock one is given by the following Black-Scholes type model:
\begin{equation}\label{bs}
\dx S_t = S_t[\;\mu \dx t + \sigma_1\dx W_1(t) + \sigma_2\dx W_2(t)\;].
\end{equation}
Here the drift $\mu$ is a real constant and the volatilities $\sigma_1, \sigma_2$ are any two positive numbers and $W_1$ and $W_2$ are independent Brownian motions. The stochastic differential equation driving $S^{'}$ is left unspecified. We only assume that it is a strong solution of a differential equation involving $W_1$ and $W_2$.
To generate incompleteness, we assume that trading is allowed only in stock one and not on stock two. 

Now suppose we want to hedge a contingent claim $C$ by trading in stock one. If we start with an initial investment of $x$ and follow a trading strategy $\pi$, the wealth at the end of the trading peiod is given by 
\[
W_{\Tau}(x,\pi) = x + \int_0^{\Tau} \pi_t\dx S_t.
\]
The quantity $(C - W_{\Tau}(x,\pi))^+$ is known as shortfall. In superhedging, we guarantee to have a shortfall of zero almost surely. This, however, needs a large initial amount $x$ which sometimes investors are unable to meet. Thus it makes sense to allow shortfall in such a way that the risk is not too large.

One common way is to fix a small number $\alpha$ as the level of endurance and allow such strategies such that the $q$th. moment of the  shortfall is bounded above by $\alpha$. That is to say,
\begin{equation}\label{momentconstraint}
\e\left[(C - W_{\Tau}(x,\pi))^+\right]^q \leq \alpha
\end{equation}
for some $q \geq 1$. Our objective is then to find the minumum real $x$ which allows us to satisfy~(\ref{momentconstraint}). 

Such a problem can be easily formulated as in subsection~\ref{statement} by a suitable choice of convex risk measure. This has been done in detail in \cite{follmer}, pages 212-218, where the reader can look for the proofs.

The sample space may be any probability space $\Omega$ on which a two dimensional Brownian motion is defined. The filtration is the augmented Brownian filtration and $P$ is the Wiener measure on this filtered probability space. We consider 
\begin{equation}\label{lexmp2}
\cxl = \{ X \in\ltwo\;\Big\vert\; P(X\ge 0)=1,\;\;\e(X)=1 \}
\end{equation}
and for $X \in \cxl$, define
\begin{equation}\label{fexmp2}
\cxf(X) := \e(XC) - {(q\alpha)}^{1/q}\norm{X}_p. 
\end{equation}
Here $p$ is given by ${1}/{p} + {1}/{q} = 1$.
We can only solve the problem for a finite $q$ greater than $2$. For such a $q$, it is immediate that $p \in (1,2)$. With this definition, determining the price of the option is the same problem as stated in equation~(\ref{*}).

\noindent{\sc {remark}.} We have taken $\cxl$ in~(\ref{lexmp2}) to be a subset of $\ltwo$ which is not usual (see \cite{follmer}, pages 212-218). However, as long as $C$ has more than two moments, this does not create any additional troubles.
\smallskip

First, we need to determine the elements of $\cinth$ defined by~(\ref{cinthperp}). Since trading is allowed only on stock one, it suffices to find the probability measures in $\cxl$ under which $S$ is a martingale in $[0,\Tau]$. The standard tool for such problems is to use  Girsanov's Theorem. Let $Q$ be a measure equivalent to $P$ under which $S$ is a martingale. Without loss of generality, we can assume that  
\begin{equation}\label{expmgle}
N_t=\e\left[ \dqdp \;\Big\vert \mathcal{F}_t\right]=\exp( L_t - 1/2 \langle L \rangle_t ),
\end{equation}
for some $L$ which is a local martingale and $\langle.\rangle$ refers to the quadratic variation of $L$. Then, by Girsanov's Theorem, if $M_t$ is a martingale under the original measure $P$, the process $\overline{M}$, given by
\begin{equation}\label{girsanov}
\overline{M}:= M - \langle M, L\rangle,
\end{equation}
is a local martingale under the new measure $Q$. Here $\langle M, L\rangle$ refers to the mutual variation between the two processes $M$ and $L$.

Now the process $(W_1,W_2)$ is a two dimensional Brownian motion. Construct a new pair of independent Brownian motions by the following rotation:
\[
\widetilde{W}_1 = \frac{\sigma_1 W_1 + \sigma_2 W_2}{\sqrt{\sigma_1^2 + \sigma_2^2}},\;\; \; 
\widetilde{W}_2 = \frac{-\sigma_2 W_1 + \sigma_1 W_2}{\sqrt{\sigma_1^2 + \sigma_2^2}}
\]
Clearly, $(\widetilde{W}_1, \widetilde{W}_2)$ is another two dimensional Brownian motion which generates the same filtration as $(W_1, W_2)$. By the {\it Predictable Representation Property} of the Brownian filtration, one can write the local martingale $L$ in equation~(\ref{expmgle}) as
\begin{equation}\label{whatisl}
\dx L_t = z_t\dx \widetilde{W}_1(t) + y_t \dx\widetilde{W}_2(t),
\end{equation}
for some progressively measurable processes $z$ and $w$. Thus the martingale $N$, in~(\ref{expmgle}), can be written in another form
\begin{equation}\label{decompofmes}
\dx N_t = N_t[\; z_t\dx \widetilde{W}_1(t) + y_t \dx\widetilde{W}_2(t)\;] = \dx N_1(t) + \dx N_2(t), 
\end{equation}
where $N_1$ and $N_2$ are local martingales with $\langle N_1, N_2 \rangle \equiv 0$. Now, by equation~(\ref{girsanov}), under the new measure $Q$, the process given by
\[
\dx\overline{W}_1(t) := \dx\tw_1(t) - z_t \dx t
\]
is a new Brownian motion. One can write the stochastic differential equation for $S$, as in equation~(\ref{bs}), in terms of $\ow_1$ in the following way
\[
\dx S_t = S_t[\; (\mu + \sigma^*z_t)\dx t + \sigma^* \dx\ow_1(t)\;]  
\]
where $\sigma^* = \sqrt{\sigma_1^2 + \sigma_2^2}$.
Thus, if under the new measure $Q$, the process $S$ is a martinagle, the only solution of $z$ is given by
\begin{equation}\label{solnofz}
z_t \equiv -\frac{\mu}{\sigma^*}.
\end{equation}
Since $S$ is adapted to the filtration of $\tw_1$ alone, the process $y_t$ can be any progressively measurable process which makes $\int y\dx\tw_2$ a true martingale. Once the processes $z$ and $y$ are described, through equations~(\ref{whatisl}) and~(\ref{expmgle}), we have characterised the class $\cinth$ of all the martingale measures for $S$. 

We are now ready to solve the problems of hedging. Specifically, we  take the example of the following European options with strike $M$, whose returns at the terminal time is
\begin{equation}\label{options}
C = (S_{\Tau} - M)^+.
\end{equation}
As discussed before, we consider $\cxl$ and $\cxf$ as given by~(\ref{lexmp2}) and~(\ref{fexmp2}), and try to use Theorem~\ref{finalthm}. We still meet some difficulties: $\cxl$ is not bounded in norm, as required by  Theorem~\ref{finalthm}. However, we can truncate or localise the problem in the following way. For any large $k$ let $\ltk := \{ X\in \ltwo,\;\;\norm{X}_2 \leq k \}$, and define
\begin{equation}
\cxl_k := \cxl \cap \ltk,\;\;\; \cxf_k(Q) := \cxf(Q),\;\; Q\in\cxl_k.        
\end{equation}
Here $\cxl$ is defined in~(\ref{lexmp2}) and $\cxf$ is defined in~(\ref{fexmp2}) with $C$ as in~(\ref{options}).
This also makes a certain intuitive sense as a penalty corresponding to a risk measure, see~(\ref{rhox}), since we are putting a heavy penalty of $\infty$ to measures $Q$ which are far away from $P$ in the sense that $\norm{{\dx Q}/{\dx P}}_2 > k$.
Since we are interested in large values of the parameter $k$, we can assume that the set $\cinth$ is non-empty. Also, since the random variable $(S_{\Tau} - M)^+$ has all moments finite, the functional $\cxf_k$ is clearly lipschitz with respect to the $\lp$ norm, satisfying assumption~(\ref{lplip}) of Theorem~\ref{finalthm}. The filtration is the augmented Brownian filtration generated by the two-dimensional Brownian motion $(W_1, W_2)$. Thus, all martingales with respect to this filtration have continuous versions. A direct application of Theorem~\ref{finalthm} would give us the following result.
\smallskip

\noindent{\sc result.} For any $y\in\mathbb{R}$, it is possible to have a self financing trading strategy $\pi$ such that
\[
\e^Q\left(y + \int_0^{\Tau} \pi_u \dx S_u\right) \geq \cxf(Q),
\]
for all $Q\in\cxl_k$, if and only if
\begin{equation}\label{pricek}
y \geq \sup_{Q\in\cinth\cap B_k}\cxf(Q).
\end{equation}

Obviously, as $k$ tends to infinity, $\cxl_k$ and $\cxf_k$ tends to $\cxl$ and $\cxf$ respectively. The value on the right-hand-side of~(\ref{pricek}) thus increases to
\[
\sup_{Q\in\cinth} \left[ \e^Q(S_{\Tau} - M)^+ - {(q\alpha)}^{{1}/{q}}\norm{{\dx Q}/{\dx P}}_p\right].
\]
We can define this limiting value to be the sellers' price of the option, since this is the infimum amount from which one can start and approximately hedge the contingent claim in the sense of~(\ref{momentconstraint}). However, from our proofs it is not apparent if there is a strategy which achieves it. 

\section{Conclusion} We consider the problem of attaining acceptability by trading under convex constraints. We start with an arbitrary convex collection of scenario measures and corresponding floors, and determine the minimum capital required so that the terminal wealth can be made acceptable. Our main result states that the minimum capital is equal to the supremum of the floors over all such scenarios under which the stock price process is a martingale. We show in an example how such a result can determine the capital requirement for hedging a contingent claim with controlled shortfall.

\section{Appendix}

\noindent{\sc{proof of lemma~\ref{Aisbdd}.}} We shall first show that for any $\epsilon > 0$, there is a $\delta > 0$ such that for any $X\in \cxl$,
\begin{equation}\label{continuity}
\norm{T(X)} < \delta \;\Rightarrow\; d_p(X,\cinth) < \epsilon.
\end{equation}
We shall prove this by contradiction. Fix an $\epsilon > 0$, suppose that~(\ref{continuity}) does not hold for any $\delta > 0$. Thus for every $\delta_n = 1/n$, one can find $X_n \in \cxl$ such that $\norm{T(X_n)} < \frac{1}{n}, \;\;\mathrm{but}\;\; d_p(X_n,\cinth)\geq \epsilon$. But this is clearly impossible by~(\ref{pcont}).

Now take any $X\in\cxl$. Since $\Pi$ is the $\lp$ closure of $\cinth$, it is clear that
\begin{equation}\label{distclos}
d_p(X,\cinth) =  d_p(X,\Pi) = \norm{X - S_{\Pi}(X)}_p.
\end{equation}
From~(\ref{scaling}) we know that for any $\alpha \geq 0$, if we denote 
\begin{equation}\label{xalpha}
X_{\alpha} := \alpha X + (1-\alpha) S_{\Pi}(X),
\end{equation}
we have $S_{\Pi}\left(X_{\alpha}\right)= S_{\Pi}(X)$. Thus
\begin{eqnarray}\label{dpscaling}
d_p(X_{\alpha}, \cinth) = \norm{X_{\alpha} - S_{\Pi}(X_{\alpha})}_p 
&=& \norm{\alpha X + (1-\alpha) S_{\Pi}(X) - S_{\Pi}(X)}_p\nonumber\\
= \norm{\alpha(X-S_{\Pi}(X))}_p &=& \alpha \norm{X-S_{\Pi}(X)}_p = \alpha. d_p(X,\cinth).
\end{eqnarray}

Since $S_{\Pi}(X)$ is an element of $\Pi$, which is the $\lp$ closure of $\cinth$, we can choose a sequence of elements $Y_n \in \cinth$ such that $\norm{Y_n - S_{\Pi}(X)}_p \rightarrow 0$.
For an $\alpha < \delta/\norm{T(X)}$ we would have
\[
\norm{T\left(\alpha X + (1-\alpha)Y_n\right)} = \alpha\norm{T(X)} < \delta.
\]
Thus, from condition~(\ref{continuity}) we get that
\begin{equation}\label{stray1}
d_p(\alpha X + (1-\alpha)Y_n, \cinth) < \epsilon,\;\forall\;n\in\mathbb{N}.
\end{equation}
However, by the triangle inequality, we have
\begin{eqnarray*}
d_p(X_{\alpha}, \cinth) &\leq& \limsup_{n\rightarrow\infty}\left[\norm{X_{\alpha} -\alpha X - (1-\alpha)Y_n}_p + d_p(\alpha X + (1-\alpha)Y_n, \cinth)\right]\\
&\le& (1-\alpha)\limsup_n\norm{S_{\Pi}(X) - Y_n}_p + \limsup_nd_p(\alpha X + (1-\alpha)Y_n, \cinth)\\
&=& 0 + \limsup d_p(\alpha X + (1-\alpha)Y_n, \cinth) \leq \epsilon.
\end{eqnarray*}
The last inequality is due to~(\ref{stray1}). Thus for any $\alpha < \delta/\norm{T(X)}$, we have
\begin{equation}\label{dpxalphaeps}
d_p(X_{\alpha},\cinth) \leq \epsilon.
\end{equation}

Now we prove~(\ref{plipt}). If $X\in \cxl$ is such that $\norm{T(X)}=0$ then $X \in G^{\perp}$ and thus $X \in \cinth$. Hence $X=S_{\Pi}(X)$ and $\norm{X-S_{\Pi}(X)}_p= 0$ and ~(\ref{plipt}) is obviously satisfied. If $\norm{T(X)}\neq 0$, choose $\alpha=\delta/(2 \norm{T(X)})$. Applying ~(\ref{dpxalphaeps}), we infer $d_p(X_{\alpha},\cinth) \leq \epsilon$. By taking $M_1=2\epsilon/\delta$, we see that~(\ref{dpscaling}) implies
\[
d_p(X,\cinth)=\frac{1}{\alpha}d_p(X_{\alpha},\cinth) \leq \frac{\epsilon}{\alpha}= M_1\norm{T(X)}.
\]
This proves ~(\ref{plipt}) and hence the Lemma.\hfill $\Box$
\medskip

\noindent{\sc{proof of lemma~\ref{Bisbdd}.}} If $X \in \mcal{Z}$, then $T(X)=0$ and $S_{\Pi}(X)=X$. Thus if~(\ref{blipt}) holds for some $z\in\mathbb{R}$, we must have $z \geq \cxf(X)$. Taking supremum over all $X \in \cinth$, we infer $z \geq \sup_{X\in \cinth} \cxf(X)$.

Sufficiency follows, since for any $z \geq \sup_{X\in \cinth} \cxf(X)=\sup_{X\in \Pi} \cxf(X)$, the left-hand side of ~(\ref{blipt}) is non-positive, and we can take $M_2=0$. \hfill $\Box$
\smallskip

\noindent{\sc proof of lemma~\ref{approxmart}.} Let us remember that $T$ is the projection operator onto the subspace $\overline{G}$. Thus $T(Y_n) \rightarrow 0$ implies that there is a sequence $\{Z_n\}\subseteq G^{\perp}$, such that
\begin{equation}\label{cvg1}
\lim_{n\rightarrow \infty}\norm{Y_n - Z_n} = 0.
\end{equation}
Hence, it also follows that
\begin{equation}\label{expconv}
\lim\e(Z_n) = \lim\e(Y_n) =1.
\end{equation}
Recall that $\e(Y_n)=1$ for all $n$, simply by virtue of being a member of $\cxl$. Thus, if we define $c_n := \e(Z_n)$ then, by~(\ref{expconv}), $c_n \rightarrow 1$, and hence is non-zero for all $n > N_2$, for some $N_2\in\mathbb{N}$. Thus, for all $n>N_2$, we can define $M_n:= c_n^{-1}Z_n$ to get 
\begin{equation}\label{whatismn}
\e(M_n)=1.
\end{equation}
Now, since $\sup_n\norm{Y_n} \leq \mcal{K}$ by assumption~(\ref{describecxl}), and~(\ref{cvg1}) holds, the sequence $\{Z_n\}$ is also uniformly bounded in the $\ltwo$ norm. Hence, it follows that
\begin{equation}\label{cvg2}
\norm{Y_n - M_n}\leq \norm{Y_n - Z_n} + \left(1-c_n^{-1}\right)\norm{Z_n} \rightarrow 0.
\end{equation}
As a corollary of the limit in~(\ref{cvg2}), we infer that given any $\epsilon > 0$, there is a $N_3$ such that 
\begin{equation}\label{eps2ball}
\norm{M_n} \leq \norm{Y_n} + \norm{Y_n-M_n} < \mcal{K} + \epsilon, \;\; \forall \;n > N_3.
\end{equation}

Now $\{M_n\} \subseteq G^{\perp}$ implies $T(M_n)=0$. Thus if we define
\[
M_n(t) := \e\left [ M_n \vert \mcal{F}_t \right ],
\]
then, by Lemma~\ref{zeroproj}, the process
\begin{equation}\label{whatisy}
Y_n(t):=S(t).M_n(t)
\end{equation}
is a martingale under $P$ in the time interval $[0,\Tau]$. Note that by assumption 2 of Theorem~\ref{finalthm}, we can choose a continuous version of $M_n(t)$. Since we assume $\mcal{F}_{\tau}$ to be the entire $\sigma$-algebra, we identify
\begin{equation}\label{toteeornottotee}
M_n(\tau) = M_n.
\end{equation}
Also, by our normalisation in~(\ref{whatismn}), we note that
\begin{equation}\label{startone}
M_n(0) = \e(M_n)=1.
\end{equation}
Let $\sigma_n$ be the stopping time defined by
\begin{equation}\label{sigman}
\sigma_n := \inf\left\{t \;\vert\; M_n(t) = 0\right \} \wedge \Tau.
\end{equation}

\noindent{\bf Claim.} We shall defer the proof of the following claim:
\begin{equation}\label{cvg3}
\norm{ M_n - M_n({\sigma_n})}_p \rightarrow 0\;\;\mathrm{as}\; n\rightarrow \infty.
\end{equation}
Assuming that the above claim is true, note that, since $p > 1$,~(\ref{cvg3}) implies
\begin{equation}\label{expecconvg}
\lim\e\left(M_n(\sigma_n)\right) = \lim\e(M_n) = 1.
\end{equation}
Thus, as before, there exists $N_4 \in \mathbb{N}$, such that for all $n>N_4$, if we define $d_n \stackrel{\triangle}{=} \e(M_n(\sigma_n))$,
then the following random variables are well-defined
\begin{equation}\label{lnscale}
L_n \stackrel{\triangle}{=} d_n^{-1}M_n(\sigma_n),\;\;\; \e(L_n) = 1.
\end{equation}
Since the martingale $M_n(t)$ is continuous, by~(\ref{startone}) and our choice of $\sigma_n$ in~(\ref{sigman}), we see that
\begin{equation}\label{mnnonneg}
M_n(\sigma_n) \geq 0, \;\mathrm{a.s.}\;P,
\end{equation}
and
\begin{equation}\label{stopyn}
Y_n({t\wedge\sigma_n}) = M_n({t\wedge\sigma_n})S(t),\;\forall\; 0\leq t\leq\tau.
\end{equation}
However, by the optional sampling theorem, the process on the left-hand side of the above expression is an $\mcal{F}_t$-martingale. Thus, 
the process of the right-hand side of~(\ref{stopyn}) is also an $\mcal{F}_t$-martingale. For every $n > N_4$, note that $d_n^{-1}M_n({t\wedge\sigma_n}) = \e\left(L_n\;\vert\;\mcal{F}_t\right)$, and hence
\begin{equation}\label{scmgle}
\{ \e\left(L_n\;\Big\vert\;\mcal{F}_t\right).S(t),\mcal{F}_t\}_{0\leq t \leq \tau}
\end{equation}
is also a martingale. 

By~(\ref{lnscale}) and~(\ref{mnnonneg}), we can change the measure $P$, by defining
\[
\dx Q_n/\dx P \stackrel{\triangle}{=} L_n,\;\;\forall\;n>N_4.
\]
Then from~(\ref{scmgle}) one can use Bayes' rule in the reverse direction to conclude that $Q_n$ is a sequence of martingale measures. Or, in other words, from Lemma~\ref{zeroismgle}, we conclude that
\[
L_n \in G^{\perp}, \;\forall\; n > N_4 .
\]
To prove Lemma~(\ref{approxmart}), now we only need to show that conditions~(\ref{stepone1}) and~(\ref{stepone2}) hold. The process $\{M^2_n(t),\mcal{F}_t\}$ is a submartingale for every $n$, and hence we have
\begin{equation}\label{doob}
\norm{M_n({\sigma_n})} \leq \norm{M_n}.
\end{equation}
Also from~(\ref{expecconvg}), it is immediate that $d_n \rightarrow 1$ and hence
\begin{eqnarray*}
\limsup_n\norm{L_n} &\leq& \lim d_n^{-1}.\limsup_n\norm{M_n({\sigma_n})}\\
&\leq& \limsup_n\norm{M_n}= \limsup_n\norm{Y_n} \leq \mcal{K}.
\end{eqnarray*}
The only equality above is due to~(\ref{cvg2}) and the final inequality is from~(\ref{describecxl}). This clearly proves condition~(\ref{stepone2}). To prove, condition~(\ref{stepone1}), notice that, by the triangle inequality, $\lim_{n\rightarrow\infty}\norm{Y_n-L_n}_p$ is bounded above by
\begin{equation}\label{gotozero}
{\limsup}\norm{Y_n-M_n}_p + {\limsup}\norm{M_n-M_n(\sigma_n)}_p
 + {\limsup}\norm{M_n(\sigma_n)-L_n}_p
\end{equation}
The first term is zero by~(\ref{cvg2}). The second term is zero by~(\ref{cvg3}). For the third term, an application of~(\ref{holderyield}) and~(\ref{doob}) will show that it is less than
\begin{eqnarray*}
\limsup\norm{M_n(\sigma_n)-L_n} &\leq& \limsup\left[\left(1-d_n^{-1}\right)\norm{M_n(\sigma_n)}\right]\\
\leq \limsup\left[\left(1-d_n^{-1}\right)\norm{M_n}\right]
&=& (\mcal{K}+\epsilon)\limsup\left(1-d_n^{-1}\right) =0.
\end{eqnarray*}
The limiting bound on $\norm{M_n}$ is obtained from~(\ref{eps2ball}). This proves that the left-hand side of~(\ref{gotozero}) is zero. We have thus shown condition~(\ref{stepone1}) holds and hence Lemma~\ref{approxmart} is proved. \hfill $\Box$
\medskip

\noindent{\sc proof of claim~(\ref{cvg3}).} Finally it remains to prove~(\ref{cvg3}). Note that by continuity of the martingale $M_n(t)$, we have $M_n({\sigma_n})=0$ on the set $\{\sigma_n < \Tau\}$. Also, due to~(\ref{toteeornottotee}), on the event $\{\sigma_n = \Tau\}$, both the random variables $M_n$ and $M_n(\sigma_n)$ are the same. Combining, we get
\begin{equation}\label{aaa}
\e | M_n - M_n({\sigma_n}) |^p = \e\left[ | M_n - M_n({\sigma_n})|^p 1_{\{\sigma_n < \Tau\}}\right] = \e\left[ | M_n|^p 1_{\{\sigma_n <\Tau \}}\right]
\end{equation}

Fix an $\epsilon > 0$. The last term above can be expressed as:
\begin{eqnarray}\label{3terms}
\e\left[ |M_n|^p 1_{\{\sigma_n < \Tau \}}\right] &=& 
\e\left[ |M_n|^p 1_{\{\sigma_n < \Tau \} \cap \{M_n > \epsilon\}}\right] \nonumber \\
&+& \e\left[ | M_n|^p 1_{\{\sigma_n < \Tau \}\cap \{M_n < -\epsilon\}} \right]\nonumber \\
 &+& \e\left[ |M_n|^p 1_{\{\sigma_n < \Tau \}\cap \{\vert M_n\vert < \epsilon\} }\right].
\end{eqnarray} 
The final term on the right-hand side of the above equation is bounded as $\e\left[ | M_n|^p 1_{\{\sigma_n < T \}\cap \{\vert M_n\vert < \epsilon\} }\right] \leq \epsilon^p$.
The second term on the right-hand side of~(\ref{3terms}) can be bounded as
\begin{equation}\label{term2bnd}
\e^P\left[ | M_n|^p 1_{\{\sigma_n < T \}\cap \{M_n < -\epsilon\}} \right] \leq \e^P\left[ |M_n|^p 1_{\{M_n < -\epsilon\}} \right].
\end{equation}
Now, by assumption in Lemma~\ref{approxmart}, the sequence $\{Y_n\}$ is a sequence in $\cxl$. Hence, by~(\ref{describecxl}), we have $P(Y_n \geq 0) =1$. Thus, on the set $\{M_n \leq 0\}$, we have $\vert M_n \vert \leq \vert Y_n - M_n \vert, \;\mathrm{a.s.}\; P$. 
The right-hand side of~(\ref{term2bnd}) can then be bounded above by 
\begin{eqnarray*}
\e^P\left[ |M_n|^p 1_{\{M_n < -\epsilon\}} \right] &\leq& \e^P\left[ |Y_n - M_n|^p 1_{\{M_n < -\epsilon\}} \right]\\
&\leq&\e^P|Y_n - M_n|^p = \left(\norm{Y_n - M_n}_p\right)^p\\
&\leq&\norm{Y_n - M_n}^p, \;\mathrm{by ~(\ref{holderyield})},
\end{eqnarray*}
which goes to zero by~(\ref{cvg2}). In the next paragraph, we shall show that the first term on the right-hand side of~(\ref{3terms}) goes to zero. Thus, combining limits of all three terms in~(\ref{3terms}), and using~(\ref{aaa}), we get that $\limsup_{n\rightarrow \infty} \e^P\left \vert M_n - M_n({\sigma_n}) \right \vert^p \leq \epsilon^p$.
Since the inequality above  holds for all $\epsilon > 0$, we have proved~(\ref{cvg3}).
\smallskip

Finally we shall show that the first term on the right-hand side of~(\ref{3terms}) goes to zero i.e., 
\begin{equation}\label{ccc1}
\lim_{n\rightarrow\infty}\e(|M_n|^p1_{\{\sigma_n < T \}\cap\{ M_n > \epsilon \}})= 0.
\end{equation}
For $r=2/p$, by~(\ref{eps2ball}), we get $\sup_n\e(|M_n|^p)^r = \sup_n\e(|M_n|^2) = \sup_n(\norm{M_n})^2$ is finite.
Since $r > 1$ by choice of $p$ ($p < 2$), this shows that the random variables $\{|M_n|^p\}_{n\in\mathbb{N}}$ is uniformly integrable. Observe that the non-negative random variables
\[
D_n \stackrel{\triangle}{=} |M_n|^p1_{\{\sigma_n < T \}\cap\{ M_n > \epsilon \}}
\]
clearly satisfy $D_n \leq |M_n|^p$, for all $n\in\mathbb{N}$. Thus the collection of random variable $\{D_n\}_{n\in\mathbb{N}}$ are also uniformly integrable. Hence, to prove~(\ref{ccc1}), it suffices to show 
\begin{equation}\label{ccc}
\lim_{n\rightarrow\infty}P({\{\sigma_n < T \}\cap\{ M_n > \epsilon \}})= 0.
\end{equation}

We shall prove~(\ref{ccc}) by contradiction. So, let us suppose that~(\ref{ccc}) does not hold, i.e., there is a $\delta > 0$ such that for a subsequence $\{n_k\}\subseteq\mathbb{N}$ we have
\begin{equation}
P(\{\sigma_{n_k} < T \}\cap\{ M_{n_k} > \epsilon \}) > \delta, \;\; \forall \;k\in\mathbb{N}. 
\end{equation}
To keep notations simple, let us do away with the subsequence notation $\{n_k\}$ and assume instead
\begin{equation}\label{ddd}
P(\{\sigma_n < T \}\cap\{ M_n > \epsilon \}) > \delta, \;\; \forall n \in \mathbb{N}. 
\end{equation}
On the event $\{ \sigma_n < T\}$, by the Optional Sampling Theorem, we have
\[
\e\left ( M_n \vert \mcal{F}_{\sigma_n} \right ) = M_n({\sigma_n}) = 0,\;\;\mathrm{a.s.}\;\; P.
\]
Thus we get the following equality
\begin{equation}\label{bbb}
P\left(\{{\sigma_n < T \}\cap\{ M_n > \epsilon \}}\right) \le P\left({\{\e\left ( M_n\vert \mcal{F}_{\sigma_n} \right )=0\}} \cap{\{ M_n > \epsilon\}} \right).
\end{equation}
Define the following non-negative random variables
\begin{equation}\label{xn}
I_n := 1_{\{\e\left( M_n \vert \mcal{F}_{\sigma_n} \right )=0\}},\;\;
J_n := P( M_n > \epsilon \vert \mcal{F}_{\sigma_n}),\;\;
K_n := I_nJ_n.
\end{equation}
Note that by conditioning the event $\{M_n > \epsilon\}$ on $\mcal{F}_{\sigma_n}$, we get
\[
P\left({\{\e\left ( M_n\vert \mcal{F}_{\sigma_n} \right )=0\}} \cap{\{ M_n > \epsilon\}} \right) = \e(I_nJ_n)=\e(K_n), 
\]
forall $n\in\mathbb{N}$. Then, by~(\ref{bbb}) and assumption~(\ref{ddd}), we have
\begin{equation}\label{expecdelta}
\e(K_n) > \delta, \;\; \forall n\in\mathbb{N}.
\end{equation}
Note that $K_n$ is a non-negative random variable, and one can get a lower bound on the tail probability by using the following basic inequality, often known as the \emph{second moment method}:
\begin{equation}
P\left(K_n \geq \frac{1}{2}\e(K_n)\right) \geq \frac{1}{4}\frac{\e(K_n)^2}{\e(K_n^2)}.
\end{equation}
Thus, combining with~(\ref{expecdelta}), we infer
\begin{equation}\label{smm}
P\left(K_n \geq \frac{\delta}{2}\right) \geq P\left(K_n \geq \frac{1}{2}\e(K_n)\right) \geq \frac{1}{4}\frac{(\e K_n)^2}{\e(K_n^2)}\geq \frac{\delta}{4}.
\end{equation}
Since $0\leq K_n\leq 1$, the last inequality follows by noting that $\e(K_n^2) \leq \e(K_n)$, and hence
\[
\frac{(\e K_n)^2}{\e(K_n^2)} \geq \e(K_n) \geq \delta.
\]
Now, note that, since $I_n$ only takes zero-one values,
\begin{equation}\label{sandvoss}
\left\{K_n \geq \frac{\delta}{2}\right\} \;\;\Leftrightarrow\;\; \left\{I_n = 1\right\}\cap\left\{J_n \geq \frac{\delta}{2}\right\}.
\end{equation}

Recall the original random variables $M_n$ which were used to define $K_n$ in~(\ref{xn}). We denote the positive and negative parts of $M_n$ by defining
\[
M_n^+ \stackrel{\triangle}{=} \max(M_n,0)\;\;\mathrm{and}\;\; M_n^-\stackrel{\triangle}{=} \max(-M_n,0).
\]
Then, on the set $\{I_n = 1\}$, we have $\e\left(M_n\vert\mcal{F}_{\sigma_n}\right)=0$, which in turn implies
\begin{equation}\label{pmbrk}
\e\left( M_n^{-} \vert \mcal{F}_{\sigma_n} \right ) = \e\left( M_n^{+} \vert \mcal{F}_{\sigma_n} \right ),\;\;\mathrm{a.s.}\;\;P. 
\end{equation}
Also, on the set $\{M_n > \epsilon\}$, we obviously have $M_n= M_n^+$, and that $\{M_n^+ > \epsilon\}$. Thus on the set $\{I_n = 1\}\cap\{J_n \geq \delta/2\}$, we have
\begin{eqnarray*}
\e\left( M_n^{-} \vert \mcal{F}_{\sigma_n} \right ) &=& \e\left( M_n^{+} \vert \mcal{F}_{\sigma_n} \right ) \geq \epsilon P(M_n^+ > \epsilon \vert \mcal{F}_{\sigma_n})\\
&=&\epsilon P(M_n > \epsilon \vert \mcal{F}_{\sigma_n}) = \epsilon J_n \;\geq\; \frac{\epsilon \delta}{2}\;\;\mathrm{a.s.}\;P.
\end{eqnarray*}
Combining the above inequality with~(\ref{sandvoss}) and~(\ref{smm}), we get that
\begin{eqnarray}\label{epsdel}
P\left( \e\left( M_n^{-} \vert \mcal{F}_{\sigma_n} \right ) \geq
\frac{\epsilon \delta}{2} \right )&\geq& P\left(\{I_n = 1\}\cap\left\{J_n \geq \frac{\delta}{2}\right\}\right)\nonumber\\
&=& P\left(K_n \geq \frac{\delta}{2}\right)\geq \frac{\delta}{4}, \;\; \forall n \in \mathbb{N}.
\end{eqnarray}

Recall the non-negative random variables $Y_n$ as in the statement of Lemma~\ref{approxmart}. Note that we always have
\begin{equation}\label{steve}
M_n^- \leq (Y_n + M_n^-)1_{\{M_n^-\neq0\}} \leq (Y_n - M_n)1_{\{M_n^-\neq0\}}\leq \vert Y_n - M_n\vert.
\end{equation}
Thus, if we let $ R_n := \e\left ( M^{n-}  \vert \mcal{F}_{\sigma_n} \right )$, from~(\ref{steve}) we conclude
\[
\e(R_n)=\e\left ( M_n^{-}\right)\leq \e\vert Y_n - M_n\vert \leq \norm{Y_n - M_n}.
\]
And thus, by ~(\ref{cvg2}), we get $\e(R_n) \rightarrow 0$. But from~(\ref{epsdel}) we get
\[
P(R_n \geq \frac{\epsilon \delta}{2}) =P\left( \e\left( M_n^{-} \vert \mcal{F}_{\sigma_n} \right ) \geq
\frac{\epsilon \delta}{2} \right ) \geq \frac{\delta}{4}, \;\; \forall n\in \mathbb{N}.
\]
This clearly contradicts $\e(R_n)\rightarrow 0$. Thus~(\ref{ddd}) cannot be true and we have thus proved~(\ref{ccc}). This completes the proof of Claim~(\ref{cvg3}).\hfill$\Box$


\begin{thebibliography}{80}

\bibitem{coherent} {\sc{Artzner, P., Delbaen, F., Eber, J.M., \& Heath, D.}} (1999)  Coherent measures of risk. \emph{Math. Finance} {\bf{9}}, 203-228.

\bibitem{barrieu} {\sc{Barrieu, P., \& El Karoui, N.}} (2005) Inf-convolution of risk measures and optimal risk transfer. \emph{Finance and Stochastics} {\bf{9}}, 269-298.

\bibitem{barrieu2} {\sc{Barrieu, P., \& El Karoui, N.}} (2005) Pricing, hedging and optimally designing derivatives via minimization of risk measures. \emph{To appear in Volume on Indifference Pricing}, Princeton University Press.

\bibitem{floor} {\sc{Carr, P., Geman, H., \& Madan, D.}} (2001)   Pricing and hedging in incomplete markets. \emph{J. Financial Economics} {\bf{62}}, 131-167.

\bibitem{karatzas} {\sc{Cvitani\'{c}, J. \& Karatzas, I.}} (1999)  On dynamic measures of risk. \emph{Finance \& Stochastics} {\bf{3}}, 451-482.

\bibitem{5authors} {\sc{Delbaen, F., Monat, P., Schachermayer, W., Schweizer, M. \& Stricker, C.}} (1997) Weighted norm inequalities and hedging in incomplete markets. \emph{Finance \& Stochastics} {\bf{1}}, 181-227. 

\bibitem{qhedge} {\sc{F\"{o}llmer, H. \& Leukert, P.}} (1999) Quantile hedging. \emph{Finance and Stochastics} {\bf{3}}, 251-273.

\bibitem{cxrisk} {\sc{F\"{o}llmer, H. \& Schied, A.}} (2002) Convex measures of risk and trading constraints. \emph{Finance and Stochastics} {\bf{6}}, 429-447.

\bibitem{follmer} {\sc{F\"{o}llmer, H. \& Schied, A.}} (2004) \emph{Stochastic Finance: An Introduction in Discrete Time}. Second Edition, Studies in Mathematics {\bf{27}}, de Gruyter, Berlin.

\bibitem{npp} {\sc{Goebel, K., \& Reich, S.}} (1984) \emph{Uniform Convexity, Hyperbolic Geometry, and Nonexpansive Mappings}. Pure and Applied Mathematics {\bf 83}, Marcel Dekker, Inc., New York.

\bibitem{karatshreve} {\sc{Karatzas, I. \& Shreve, S. E.}} (1991) \emph{Brownian Motion and Stochastic Calculus}. Second Edition, Springer-Verlag GTM {\bf 113}, New York.


\bibitem{shreve} {\sc{Larsen, K., Pirvu, T., Shreve, S., \& T\"{u}t\"{u}nc\"{u}, R.,}} (2004) Satisfying convex risk limits by trading. \emph{Finance and Stochastics} {\bf 9}, 177-195.

\bibitem{megg} {\sc{Megginson, R. E.}} (1998) \emph{An Introduction to Banach Space Theory}. Springer GTM {\bf 183}, New York.

\bibitem{protter} {\sc{Protter, P.}} (2004) \emph{Stochastic Integration and Differential Equations}. Second Edition, Stochastic Modelling and Applied Probability {\bf{21}}, Springer, pp. 244-245.

\bibitem{rhein} {\sc{Rheinl\"{a}nder, T. \& Schweizer, M.}} (1997) On {\bf{L}}$^2$-projections on a space of stochastic integrals. \emph{Annals of Probability} {\bf{25}}, 1810-1831.

\bibitem{rock} {\sc{Rockafellar, R.}} (1997) \emph{Convex Analysis}. Tenth printing and first paperback printing in the Princeton Landmarks in Mathematics series, Princeton University Press.

\bibitem{rudin} {\sc{Rudin, W.}} (1991) \emph{Functional Analysis}. Second Edition, International Series in Pure and Applied Mathematics, McGraw-Hill.

\end{thebibliography}
\end{document}